\documentstyle[fullpage, 12pt]{article}

\def\C{{\rm C \kern-.48em\vrule width.06em height.57em depth-.02em \kern.48em}}
\def\sopenC{{\rm C\kern-.15cm\vrule width.6pt height 4.1pt depth-.3pt
 \kern.15cm}}
\def\Z{{{\rm Z}\kern-.28em{\rm Z}}}
\def\N{{{\rm I}\kern-.16em{\rm N}}}
\def\R{{{\rm I}\kern-.16em{\rm R}}}

\def\set#1{\langle #1 \rangle}

\def\belowrightarrow#1{{{{}\over\ #1\ }\kern-1.1em\to}}

\def\goback#1{\setbox0\hbox{#1}\kern-\wd0 \relax}
\def\eqbd{\mathop{{:}{=}}}

\begin{document}
\large

\title{Open problems on GKK $\tau$-matrices}
\author{Olga Holtz  \\
\small Department of Computer Sciences \\
\small University of Wisconsin \\
\small Madison, Wisconsin 53706 U.S.A. \\
\small holtz@cs.wisc.edu
 \and
Hans Schneider  \\
\small Department of Mathematics \\
\small University of Wisconsin \\
\small Madison, Wisconsin 53706 U.S.A. \\
\small hans@math.wisc.edu}

\date{ 20 August 2001}
\maketitle

\abstract{We propose several open problems on GKK $\tau$-matrices
raised by examples showing that some such matrices are unstable.}

\section{Motivation}

A complex square matrix is called `(positive) stable' if all its
eigenvalues have positive real part. Such stability of a matrix,
or more usually of its negative, is a basic concern in several
fields of mathematics, e.g., in the study of ODEs. Three important
classes of matrices, namely the nonsingular totally nonnegative
matrices, the Hermitian positive definite matrices, and the
M-matrices are known to be stable. Matrices in these three classes
also share other properties: they are P-matrices, they are weakly
sign symmetric, and they satisfy eigenvalue monotonicity (for
definitions of these terms see below). In order to investigate the
relations between these properties supersets of the three basic
matrix classes defined in terms of some or all of the above
properties have been studied. Thus Gantmacher-Krein \cite{GK},
Kotelyansky \cite{K}, Fan \cite{F} and Carlson \cite{C1},
\cite{C2} investigated properties of P-matrices that are weakly
sign symmetric (the so called GKK matrices), while Engel-Schneider
\cite{ES} considered  $\tau$-matrices which are defined as
P-matrices that satisfy eigenvalue monotonicity. The latter were
further investigated by Varga \cite{V}, Hershkowitz-Berman
\cite{HB} and Mehrmann \cite{M}. The stability of GKK and $\tau$-matrices 
was conjectured and proved in low dimensions in various
papers, see Holtz \cite{OH} for details and further history and
see also the survey by Hershkowitz \cite{DH}.

The GKK $\tau$-matrices satisfy the three properties discussed above.
Nevertheless, Holtz \cite{OH} has shown that there are unstable GKK
$\tau$-matrices. This raises the question under what conditions a
GKK $\tau$-matrix is stable. It is the purpose of this note to
outline five problems related to this question.

\section{Basic notions}

Let $\set{n}$ denote the collection of all increasing sequences with elements
from the set $\{1,2,\ldots,n\}$ and let $\# \alpha$ denote the size of the
sequence $\alpha$. Given a matrix $A\in \C^{n\times n}$, we will
use the notation $A(\alpha,\beta)$ for the submatrix of $A$ whose rows are 
indexed by $\alpha$ and columns by $\beta$ ($\alpha$, $\beta\in \set{n}$) 
and $A[\alpha,\beta]$ for the minor $\det A(\alpha,\beta)$ if $\# \alpha=\#
\beta$. For simplicity, $A(\alpha)$ will stand for $A(\alpha,\alpha)$
and $A[\alpha]$ for $A[\alpha,\alpha]$. By definition,
$A[\emptyset]\eqbd 1$. Given $\alpha,\beta\in\set{n}$ with $\# \alpha=
\# \beta$, call the number $\# \alpha -\# (\alpha\cap \beta)$ the
{\em dispersal\/} of the minor $A[\alpha,\beta]$.
 A matrix $A$ is called a {\em $P$-matrix\/} if $A[\alpha]>0$
for all $\alpha\in \set{n}$.  $A$ is said to be {\em sign-symmetric\/}
if \begin{equation} A[\alpha,\beta] A[\beta,\alpha]\geq 0 \label{gkk}
\end{equation}
for all minors $A[\alpha, \beta]$. $A$ is called {\em weakly
sign-symmetric\/}~\cite{GK}  if~(\ref{gkk}) holds for all minors
$A[\alpha,\beta]$ with dispersal $1$, and these are also referred to
as {\em almost principal\/}. Weakly sign-symmetric $P$-matrices
were called {$\!$\em GKK\/} in~\cite{F} after Gantmacher, Krein,
and Kotelyansky. A $P$-matrix is {\em strict GKK\/} if the
inequalities in~(\ref{gkk}) are strict for all almost principal
minors. Let $\sigma(A)$ denote the spectrum of $A$ and  let
$$l(A)\eqbd \min \sigma(A) \cap \R,$$ with the understanding that,
in this setting, $\min \emptyset=\infty$. A matrix $A$ is called
an {\em $\omega$-matrix\/}~\cite{ES} if it has eigenvalue
monotonicity in the  sense that $$ l(A(\alpha,\alpha))\leq
l(A(\beta,\beta))<\infty \qquad \mbox{whenever} \quad \emptyset
\neq \beta \subseteq \alpha \in \set{n}.$$ $A$ is a {\em
$\tau$-matrix\/} if, in addition, $l(A) \geq 0$. A matrix is
called {\em positive  stable\/} if its spectrum lies entirely in
the open right half plane. In the sequel, we will shorten the term
`positive stable' to simply `stable'.

\section{GKK $\tau$-matrix stability and related problems}

\begin{description}
\item{\bf 1. Strict GKK matrices: closure of the set.}
What GKK matrices can be approximated arbitrarily well by strict GKK matrices?
In particular, can the matrices constructed
in~\cite{OH} be approximated by $\tau$-matrices that are strict GKK?
The matrices in~\cite{OH} themselves are not strict GKK. The negative
answer to this question would give rise to Problem~1a.
\begin{description}
\item{\bf 1a. Strict GKK matrices: stability.} Are strict GKK matrices stable? Are strict GKK  $\tau$-matrices stable?
\end{description}
\item{\bf 2. Dispersal condition sufficient for stability.}  The counterexample
given in~\cite{OH} shows that it is not
sufficient for stability of a $P$-matrix $A$ that  the inequalities~(\ref{gkk})
hold for all minors of dispersal less than or equal to $d=1$.
Carlson\index{Carlson's theorem}'s theorem~\cite{C2}  asserts that  the value
$d=n$ is sufficient for stability (in other words, that sign-symmetric
matrices are stable). What minimal value of the parameter $d$ would
guarantee stability? In particular, does that value depend on $n$?

\item{\bf 3. Classes of stable matrices.}
Following a conjecture by Varga \cite{V}, Mehrmann~\cite{M} showed
that $\tau$-matrices of order up to $4$ satisfy the inequality
 \begin{equation} \label{Varg}
|\arg(\lambda - l(A))|\leq \frac{\pi}{2}-\frac{\pi}{n} \qquad
\forall \lambda \in \sigma(A),
 \end{equation}
 which is a property stronger than stability.

 However, in \cite{OH} it was shown that the GKK $\tau$-matrices
$A_{2k+2,k,t}$ of order $2k+2$ defined there are unstable for all
integers $k,\ k
> 20$, and all sufficiently small positive $t$
(and numerical evidence suggests that there may be matrices
$A_{n,k,t}$ of smaller order that are unstable). For various
classes $\mathcal{C}$ of matrices, this raises the
following three questions concerning the subclass ${\mathcal C}_n$ 
of all matrices in $\mathcal{C}$ of order less than or equal to $n$.
\begin{description}
\item[3a.] What is the maximum $n$ such that all matrices in class
${\mathcal C}_n$ are stable?
\item[3b.]  What is the maximum $n$ such that all matrices in class
${\mathcal C}_n$ satisfy (\ref{Varg})?
\item[3c.]  For given $n$,
do the stable matrices in class ${\mathcal C}_n$ satisfy
(\ref{Varg})?
\end{description}

Specific classes $\mathcal{C}$ of interest include 
\begin{description}
\item[i.] The class of matrices $A_{n,k,t}$ for positive integers $n$ and
$k$, and for $t\in (0,1)$,
\item[ii.]
The class of GKK $\tau$-matrices,
\item[iii.]
The class of GKK matrices,
\item[iv.]
The class of $\tau$-matrices.
\end{description}

\item{\bf 4. Assignment of principal minors.} Given $n\in \N$ and
numbers
$(p_\alpha)_{\emptyset \neq \alpha\in \set{n}}$,  is there a matrix $A$
such that  $A[\alpha]=p_\alpha$ for all  $\alpha$?

This question was originally  motivated
by the Gantmacher-Krein-Carlson theorem~(\cite{GK} and \cite{C1}), which states
that a $P$-matrix $A$ is GKK if and only if its minors satisfy the generalized
Hadamard-Fischer inequality
\begin{equation}
A[\alpha] A[\beta] \geq A[{\alpha\cup \beta}]
A[{\alpha \cap \beta}] \qquad \forall \alpha,\beta \in
\set{n}. \label{HF}
\end{equation}

An answer to Problem~4, coupled with inequalities~(\ref{HF}), would
therefore
allow to decide, for a collection of positive numbers $(p_\alpha)$, whether
there exists a GKK matrix such that $A[\alpha]=p_\alpha$ for all  $\alpha$.
Since  $\sigma(A)$ is determined by the numbers $p_\alpha$, one could
then find, at least in principle, all possible spectra of GKK matrices.

This problem is in fact equivalent to a certain inverse eigenvalue problem.
Indeed, specifying all principal minors implies specifying characteristic
polynomials, and hence  all eigenvalues, of the principal submatrices, and
vice versa.

Problem~4 can be specialized to various classes of matrices, e.g.,
Hermitian or nonnegative. A part of this problem for Hermitian matrices
is solved in~\cite{BK}, where it is shown how to construct a symmetric
$m$-band matrix of order $n$ from the given (necessarily real) eigenvalues
of the $m$ leading principal submatrices of greatest order. Such a matrix
always exists (but is not unique) whenever the eigenvalues of $A(1{:}j)$
interlace those of $A(1{:}j{+}1)$, $j=n-m+1,\ldots, n-1$.  Notice that,
by  the Hermite-Biehler
theorem (see, e.g.,~\cite[p.21]{ChM}), two real polynomials $p$ and $q$
have real interlacing roots if and only if the roots of the polynomial $p+iq$ are all
on the same side of the real axis, and the latter can be checked using the
 Hurwitz matrix for the polynomial $w(z)\eqbd p(iz)+iq(iz)$, which is
constructed from the coefficients of $p$ and $q$.
Therefore, the interlacing property can be checked using only
the coefficients of characteristic polynomials of principal submatrices,
and these are sums of the given numbers $p_\alpha$.

\item{\bf 5. Newton's inequalities.}  For a matrix $A$, let
$$ c_j\eqbd\sum_{\# \alpha=j} A[{\alpha}]/{n \choose j}, \qquad
 \qquad  j=0,\ldots,n. $$
Does any GKK $\tau$-matrix $A$ satisfy the inequalities
\begin{equation} c_j^2\geq c_{j-1}c_{j+1}, \;\;\; j=1,\ldots,n-1?
 \label{newton} \end{equation}
This problem has possible subproblems, e.g.,
\begin{description} \item{\bf 5a. Newton's inequalities for
$M$-matrices.}
 \item {\bf 5b. Newton's inequalities for stable GKK
$\tau$-matrices.}\end{description}

These inequalities are known for real diagonal matrices, i.e., simply for
sequences of real numbers (see~\cite{N} and references therein),  as was
first proved by Newton. Since the numbers $c_j$ are invariant under
similarity, Newton's inequalities~(\ref{newton}) also hold for all 
diagonalizable matrices with real spectrum, and therefore also for the
closure of this set, viz. for {\em all\/} matrices with real spectrum. 

This question arose at an early stage of the work that led to~\cite{OH}
in an attempt to prove that Newton's inequalities hold for the GKK
matrices and imply their stability. When it became clear that there exist
$P$-matrices satisfying Newton's inequalities that are unstable,
the author of~\cite{OH} turned to constructing a counterexample to the
GKK matrix stability conjecture. Whether GKK matrices satisfy  Newton's
inequalities, however, is not known yet.

\end{description}

We hope that consideration of our questions will lead to further
interesting results on GKK $\tau$-matrices.

\section*{Acknowledgements} 
We thank Carl de Boor and Frank Uhlig for their critical reading of
the manuscript and for suggestions which have improved this paper.

\end{document}